# Higher Order Aitken Extrapolation with Application to Converging and Diverging Gauss-Seidel Iterations

Ababu Teklemariam Tiruneh[1]

[1]Lecturer, Department of Environmental Health Science. University of Swaziland.
P.O.Box 369, Mbabane H100, Swaziland.  Email: ababute@yahoo.com

**Abstract**

In this paper Aitken's extrapolation normally applied to convergent fixed point iteration is extended to extrapolate the solution of a divergent iteration. In addition, higher order Aitken extrapolation is introduced that enables successive decomposition of high Eigen values of the iteration matrix to enable convergence. While extrapolation of a convergent fixed point iteration using a geometric series sum is a known form of Aitken acceleration, it is shown in this paper that the same formula can be used to estimate the solution of sets of linear equations from diverging Gauss –Seidel iterations. In both convergent and divergent iterations, the ratios of differences among the consecutive values of iteration eventually form a convergent (divergent) series with a factor equal to the largest Eigen value of the iteration matrix. Higher order Aitken extrapolation is shown to eliminate the influence of dominant Eigen values of the iteration matrix in successive order until the iteration is determined by the lowest possible Eigen values. For the convergent part of the Gauss-Seidel iteration, further acceleration is made possible by coupling of the extrapolation technique with the successive over relaxation (SOR) method. Application examples from both convergent and divergent iterations have been provided. Coupling of the extrapolation with the SOR technique is also illustrated for a steady state two dimensional heat flow problem which was solved using MATLAB programming.

## 1. Introduction

Iterative solutions to systems of equations are widely employed for solving scientific problems. They have several advantages over direct methods.  For equations involving large number of unknowns iterative solutions involve operations of order ($n^2$) compared to direct solutions of order ($n^3$). In computer applications, iterative solutions require much less memory and are quite simple to program [1]. In addition iterative solutions are in many cases applicable to non-linear sets of equations.

One set of iterative methods that are in wide use is centered on generation of Krylov sub space such as the method of conjugate gradients developed by Hestens and Stiefel [2]. Such methods are guaranteed to converge in at most N steps for problems involving N unknowns. However, iterative solutions based on stationary methods such as Gauss-Seidel and others that are simpler to program and do not generate new iteration matrix are receiving greater application [3]. Iterative methods that mimic the physical processes involved such as the forward-backward method have been shown to produce solutions for some problems faster than krylov methods [4]. However, fixed point iterative methods are known to be less robust than Krylov methods and convergence is not guaranteed for ill-conditioned systems of equations.

Iterative processes involving fixed point iterations that are convergent such as Jacobi and Gauss-Seidel iterations are known to form error series that are diminishing in the form of geometric series [5,14, 16]. The geometric series sum based form of Aitken extrapolation for fixed point iteration involving the system of equations AX=B is based on the formula [6]:

$$X_\infty = X_k + \frac{E_k}{1 - \lambda}$$

Where $X_\infty$ is the Aitken extrapolation of the solution, $X_k$ is the approximation to the solution at the $k^{th}$ iteration, $E_k$ is the difference in X at the $k^{th}$ iteration ($X_{k+1} - X_k$) and $\lambda$ is the ratio of the difference in X values ( $\lambda = E_{k+1}/E_k$).

The geometric series extrapolation requires at least three points of the iteration. The extrapolation in terms of the three points at the k, k+1 and k+2 iterations takes the form [6].

$$X_\infty = \frac{X_{k+2}X_k - X_{k+1}X_{k+1}}{X_{k+2} - 2X_{k+1} + X_k}$$

Which is a known form of Aitken extrapolation.

Irons and Shrive [7] made a modification to Aitken's method for scalars (sequences) which can also be applied to individual x values of fixed point iterations such as the Gauss-Seidel iteration. The extrapolation to the fifth point uses four points. After the four point extrapolation the fixed point iteration formula is used to obtain the fifth point. In this way by joint application of both extrapolation and fixed point iteration, the method is said to be dynamic with better convergence properties. A dynamic model in such a form does not require restarting the iteration procedure as the method generates all necessary iterates from the latest extrapolation [8].

The reduced rank extrapolation method [9] extends the scalar form of Aitken extrapolation into vectors of a given dimension. The extrapolation formulation for the iteration vector is a vector parallel of Aitken's extrapolation:

$$X_\infty = X_K - (M - I)^{-1} * E_K$$

Where M is the iteration matrix of the fixed point iteration, I is the identity matrix of rank N, X are the iteration vector and $E_k$ is the vector difference in X at the $k^{th}$ iteration. The method involves computing the generalized inverse involving the second order difference vector and is time consuming for solving non-linear systems of equations of larger size.

Gianola and Schaffer [6] applied geometric series extrapolation for Jacobi and Gauss-Seidel iterations in animal models. The optimal relaxation factor was lower when solutions were extrapolated, but its value was not as critical in the case of extrapolation.

Fast Eigen vector computations that require matric inversion or decomposition are unsuitable for large size matrix problems as many of them involve operations of the order O ($n^3$). Kamvar, *et al* [10] applied Aitken extrapolation for accelerating page rank computations. They showed that Aitken acceleration computes the principal eigenvector of a Markov matrix under the assumption that the power-iteration estimate x (k) can be expressed as a linear combination of the first two eigenvectors.

Calude Breziniski and Michela Redivo Zaglia [11] proposed extension of Aitken's extrapolation into a general form involving transforming the sequence of iteration into a different form using known sequences which can lead to stabilization and convergence of the original iteration. The transformation, however, is not simple and straight forward and required further refinement.

Chebyshev acceleration is also a way of transforming the iteration sequence which, for iteration matrix of known upper and lower bound Eigen values, the transformed sequence using Chebyshev polynomials leads to convergence of the fixed point iteration. In Chebyshev acceleration, the sequence of iteration values is modified by multiplication with Chebyshev polynomials which are constructed from known or estimated ranges of Eigen values of the iteration matrix. The choice of the form of Chebyshev polynomials is such that the procedure leads to progressive reduction of the norm of the error vector through a min-max property which minimizes the maximum value that the polynomial has for the range of Eigen values specified [12]. However, Chebyshev acceleration has the drawback of the need to accurately estimate the bounds of Eigen values of the iteration matrix, because outside the domain of Eigen values the polynomial shows divergence and the min-max property does not hold.

**2. Method development**

For solving a system of linear equations using fixed point iteration, the Aitken extrapolation formula can be written in the form:

$$x_\infty = x_k + \frac{e_k}{1 - \lambda_k}$$

Where:
    $x_\infty$ = The estimate of the solution at the limit of iteration
    $e_k$ = The difference in consecutive x values, i.e., $x_{k+1} - x_k$
    $\lambda_k$ = The ratio of differences in x values, i.e.,

$$\lambda_k = \frac{e_{k+1}}{e_k} = \frac{x_{k+2} - x_{k+1}}{x_{k+1} - x_k}$$

It will now be shown that the above formula is applicable to both convergent as well as divergent Gauss-Seidel iteration. In addition it will also be shown that successive application of the Aitken extrapolation formula to a higher order will result in deflation of the dominant Eigen values one by one there by transforming a divergent iteration to a convergent form.

Let the system of linearized equations for a given problem be represented in the matrix form:

$$AX = B$$

Where A is the coefficient matrix, B is the right hand side vector and X is the solution vector. Writing the matrix A further in terms of the components L, U and D matrices gives;

$$A = U + L + D$$

Where U and L are the upper and low triangular matrices respectively and D is the diagonal matrix. For the Gauss-Seidel iteration, the system of equations now can be written as:

$$(U + L + D)X = B$$

$$(L + D)X = -UX + B$$

Using the k+1$^{th}$ and k$^{th}$ iteration X-vectors for the left and right hand sides of equation above respectively;

$$(L + D)X_{k+1} = -UX_k + B$$

Solving for $X_{k+1}$;

$$X_{k+1} = -(L + D)^{-1} * UX_k + (L + D)^{-1}B \qquad \text{............ [1]}$$

This is the Gauss-Seidel iteration and the matrix;

$$-(L + D)^{-1} * U$$

is called the iteration matrix. The actual iteration in terms of the x values (scalars) takes the form;

$$x_i^{k+1} = \frac{b_i}{a_{ii}} - \left(\frac{1}{a_{ii}}\right) \sum_{j=1}^{i-1} a_{ij} x_j^{k+1} - \left(\frac{1}{a_{ii}}\right) \sum_{j=i+1}^{n} a_{ij} x_j^{k} \qquad \text{............ [2]}$$

For the Gauss-Seidel iteration the vector of x values is successively computed from the fixed-point iteration process by:

$$X_{k+1} = -(L+D)^{-1} * UX_k + (L+D)^{-1}B = MX_k + N \quad \text{............ [3]}$$

Where the matrix $M = -(L+D)^{-1} * U$ is the iteration matrix and $N = (L+D)^{-1}B$.

The differences in the solution vector $X^k$ among consecutive steps of iteration are formulated as follows;

$$X^{k+1} = -(L+D)^{-1} * UX^k + (L+D)^{-1}B$$

$$X^{k+2} = -(L+D)^{-1} * UX^{k+1} + (L+D)^{-1}B$$

$$X^{k+2} - X^{k+1} = -(L+D)^{-1} * U(X^{k+1} - X^k)$$

In other words,

$$E_{k+1} = [-(L+D)^{-1} * U] E_k \quad \text{............ [4]}$$

For a convergent Gauss-Seidel iteration, the differences in x values $E_k$ written in terms of the difference between consecutive vectors of x values in equation (4) above converges to the zero difference vectors when the solution X vector is reached. Denoting by M the iteration matrix;

$$M = -[(L+D)^{-1} * U]$$

and the difference vector of X among consecutive steps of the iteration by;

$$[E_0], \quad [E_1], \quad [E_2], \quad \text{... ...} [E_k], \quad [E_{k+1}] \text{ ... ....}$$

$$[E_{k+1}] = M * [E_k]$$

Assuming that the initial difference in x values (here after termed the error vector) $E_0$ can be written as a linear combination of the Eigen vectors $v_1$, $v_2$, $v_N$ of the iteration matrix M of size N with corresponding Eigen values $\lambda_1$, $\lambda_2$ .......$\lambda_N$ in decreasing order of magnitude where $\lambda_1$ is the dominant Eigen value;

$$E_0 = X_1 - X_o = c_1 v_1 + c_2 v_2 + c_3 v_3 + \text{... ...} + c_N v_N$$

$$E_1 = ME_0 = c_1 M v_1 + c_2 M v_2 + c_3 M v_3 + \text{... ...} + c_N M v_N$$

$$E_1 = c_1\lambda_1 v_1 + c_2\lambda_2 v_2 + c_3\lambda_3 v_3 + \ldots\ldots + c_N\lambda_N v_N$$

$$\vdots$$

$$E_k = c_1\lambda_1^k v_1 + c_2\lambda_2^k v_2 + c_3\lambda_3^k v_3 + \ldots\ldots + c_N\lambda_N^k v_N$$

$$E_{k+1} = c_1\lambda_1^{k+1} v_1 + c_2\lambda_2^{k+1} v_2 + c_3\lambda_3^{k+1} v_3 + \ldots\ldots + c_N\lambda_N^{k+1} v_N$$

Taking the ratio of the Euclidean norms of $E_k$ and $E_{k+1}$;

$$\frac{\|E_{k+1}\|}{\|E_k\|} = \frac{\left((c_1\lambda_1^{k+1}\|v_1\|)^2 + (c_2\lambda_2^{k+1}\|v_2\|)^2 + (c_3\lambda_3^{k+1}\|v_3\|)^2 + \ldots + (c_N\lambda_N^{k+1}\|v_N\|)^2\right)^{1/2}}{\left((c_1\lambda_1^k\|v_1\|)^2 + (c_2\lambda_2^k\|v_2\|)^2 + (c_3\lambda_3^k\|v_3\|)^2 + \ldots + (c_N\lambda_N^k\|v_N\|)^2\right)^{1/2}}$$

$$= \frac{\lambda_1^{k+1}\left[(c_1\|v_1\|)^2 + \left(c_2\|v_2\|\left(\frac{\lambda_2}{\lambda_1}\right)^{k+1}\right)^2 + \left(c_3\|v_3\|\left(\frac{\lambda_3}{\lambda_1}\right)^{k+1}\right)^2 + \ldots\ldots + \left(c_N\|v_N\|\left(\frac{\lambda_N}{\lambda_1}\right)^{k+1}\right)^2\right]^{1/2}}{\lambda_1^k\left[(c_1\|v_1\|)^2 + \left(c_2\|v_2\|\left(\frac{\lambda_2}{\lambda_1}\right)^k\right)^2 + \left(c_3\|v_3\|\left(\frac{\lambda_3}{\lambda_1}\right)^k\right)^2 + \ldots\ldots + \left(c_N\|v_N\|\left(\frac{\lambda_N}{\lambda_1}\right)^k\right)^2\right]^{1/2}}$$

If $\lambda_1$ is the dominant Eigen value, then the ratios $\lambda_2/\lambda_1$, $\lambda_2/\lambda_1$, ... $\lambda_N/\lambda_1$ tend to zero at higher k values so that;

$$\frac{\|E_{k+1}\|}{\|E_k\|} = \frac{\lambda_1^{k+1}[(c_1\|v_1\|)^2\ ]^{1/2}}{\lambda_1^k[(c_1\|v_1\|)^2\ ]^{1/2}} = \lambda_1 \qquad\qquad \ldots\ldots\ldots\ldots\ [5]$$

Therefore, for both converging and diverging Gauss-Seidel iterations, the error vector ratios are determined by the dominant Eigen value of the iteration matrix, $\lambda_1$.

**Case I: Convergent Gauss-Seidel Iteration**

For a convergent iteration in which the dominant Eigen value is less than one, the difference vector $E_k$ converges to zero. Denoting the estimate of the largest Eigen value of the iteration matrix M at the $k^{th}$ iteration by $\lambda_1$;

and taking the difference among respective values of x as;

$$e_{k+1} = \lambda_1 * e_k$$

The largest Eigen value of the iteration matrix $\lambda_1$ can be expressed as;

$$\lambda_1 = \lim_{k \to \infty} \left( \frac{e_{k+1}}{e_k} \right)$$

An approximation to the dominant Eigen value $\lambda_1$ is made from the ratio of differences in x values computed at the $k^{th}$ iteration, i.e.,

$$\lambda_1 \approx \frac{e_{k+1}}{e_k}$$

Starting with the $k^{th}$ difference in x values among consecutive steps of iteration, a geometric series is formed in terms of the e values;

$$e_k, \quad \lambda_1 e_k, \quad \lambda_1^2 e_k, \quad \lambda_1^3 e_k \quad \ldots \ldots \ldots \ldots \ldots \ldots$$

The sum of this diminishing geometric series is given by (since $\lambda_1 < 1$);

$$\sum_{i=k}^{i=\infty} e_i = \frac{e_k}{1 - \lambda_1}$$

It is possible to write the e values so that;

$$e_k = x_{k+1} - x_k$$
$$e_{k+1} = x_{k+2} - x_{k+1}$$
$$e_{k+2} = x_{k+3} - x_{k+2}$$

$$\vdots$$

So that ;

$$\sum_{i=k}^{i=\infty} e_i = (x_\infty - x_k) = \frac{e_k}{1 - \lambda_1} \qquad \text{............ [6]}$$

Therefore, the solution $x_\infty$ at convergence is extrapolated from equation 6 above;

$$x_\infty = x_k + \frac{e_k}{1 - \lambda_1} \qquad \text{............ [7]}$$

**Case II: Divergent Gauss-Seidel Iteration**

For a divergent Gauss-Seidel iteration, the error propagation is studied as the iteration progresses. Let $X_0$ be the initial estimate of the solution while $X_s$ is the true solution of the system of equation. The initial error vector $E_{ro}$ can then be written as:

$$E_{ro} = X_0 - X_s$$

The difference between consecutive iteration values X is as defined before, i.e.,

$$E_k = X_{k+1} - X_k$$

The Gauss-Seidel iteration formula given in equation (3), i.e.,

$$X_{k+1} = -(L + D)^{-1} * UX_k + (L + D)^{-1}B$$

can also be rewritten as:

$$X_{k+1} = MX_k + N$$

Where $M = -(L + D)^{-1} * U$ is the iteration matrix as defined before and $N = (L + D)^{-1}B$

The successive values of X of the iteration can now be written as follows;

$$X_1 = MX_0 + N = M(X_s + E_{r0}) + N$$

$$= (MX_s + N) + ME_{ro}$$

$$= X_s + ME_{ro}$$

The above expression is true because at the true solution $X_s$, the Gauss-Seidel iteration satisfies the relation:

$$X_s = MX_s + N$$

Similarly for $X_2$;

$$X_2 = MX_1 + N = M(X_s + ME_{ro}) + N$$

$$X_2 = (MX_s + N) + M^2 E_{ro}$$

$$= X_s + M^2 E_{ro}$$

Proceeding similarly at the $k^{th}$ iteration, the X-values can be written as:

$$X_k = X_s + M^k E_{ro} \qquad \text{............ [8]}$$

In terms of the difference between consecutive x- estimates, recalling the formula derived earlier in equation (4), i.e.,

$$[E_{k+1}] = -[(L+D)^{-1} * U] * [E_k] = ME_k$$

$$X_k = X_0 + E_0 + ME_0 + M^2 E_0 + M^3 E_0 + \ldots\ldots \ldots\ldots\ldots..+M^{k-1}E_0$$

$$X_k = X_0 + \sum_{i=0}^{k-1} M^i E_0 \qquad \text{............ [9]}$$

It is seen from equations (8) and (9) above that the X values increase in proportion to the iteration matrix M.

Representing this increase in proportion to M by the largest Eigen value of the iteration matrix, $\lambda_1$, equations (8) and (9) can now be written as:

$$X_k = X_s + \lambda_1^k E_{ro} = X_s + \lambda_1^k (X_0 - X_s) \qquad \text{............ [10]}$$

$$X_k = X_0 + \sum_{i=0}^{k-1} \lambda_1^i E_0 = X_0 + \frac{E_0(\lambda_1^k - 1)}{\lambda_1 - 1} \qquad \text{............ [11]}$$

The expression in equation (11) above is derived using the general geometric series sum formula for an expanding geometric series.

Equating the expressions for $X_k$ in equations (10) and (11) above:

$$X_s + \lambda_1^k(X_0 - X_s) = X_0 + \frac{E_0(\lambda_1^k - 1)}{\lambda_1 - 1}$$

$$(\lambda_1^k - 1)(X_0 - X_s) = \frac{E_0(\lambda_1^k - 1)}{\lambda_1 - 1}$$

$$X_0 - X_s = \frac{E_0}{\lambda_1 - 1} = -\frac{E_0}{1 - \lambda_1}$$

$$X_s = X_0 + \frac{E_0}{1 - \lambda_1}$$

In general for iteration estimation made from the $k^{th}$ iteration vales of $x_k$ and individual ek values the formula can be written as:

$$x_s = x_k + \frac{e_k}{1 - \lambda_1} \qquad \text{............ [12]}$$

It can be seen, therefore, that the same formula used for extrapolating a convergent Gauss-Seidel iteration can be used to extrapolate the solution from the divergent Gauss-Seidel iteration. Such a procedure which is unconventional works well as the examples that follow illustrate.

**Higher Order Aitken Extrapolation**

For the first –order Aitken extrapolation, the ratio of the norms of the error vector was shown in equation (5) to be equal to the dominant Eigen value of the iteration matrix, i.e.,

$$\frac{\|E_{k+1}\|}{\|E_k\|} = \frac{\lambda_1^{k+1}[(c_1\|v_1\|)^2]^{1/2}}{\lambda_1^{k}[(c_1\|v_1\|)^2]^{1/2}} = \lambda_1$$

For the second order Aitken extrapolation, the extrapolation made at the $k^{th}$ and $k+1^{th}$ iterations are considered:

$$X_{s(k)}^{(1)} = X_k + \frac{E_k^{(1)}}{1-\lambda_1}$$

$$X_{s(k+1)}^{(1)} = X_{k+1} + \frac{E_{k+1}^{(1)}}{1-\lambda_1}$$

The second order error in terms of the extrapolated X vectors can be written as:

$$E_k^{(2)} = X_{s(k+1)}^{(1)} - X_{s(k)}^{(1)} = \left(X_{k+1} + \frac{E_{k+1}^{(1)}}{1-\lambda_1}\right) - \left(X_k + \frac{E_k^{(1)}}{1-\lambda_1}\right)$$

$$E_k^{(2)} = E_k^{(1)} + \frac{\Delta E_k^{(1)}}{1-\lambda_1}$$

Similarly for the $k+1^{th}$ iteration;

$$E_{k+1}^{(2)} = E_{k+1}^{(1)} + \frac{\Delta E_{k+1}^{(1)}}{1-\lambda_1}$$

At the $k^{th}$ iteration, the ratio of the norm of the error vector becomes;

$$\left(\frac{\|E_{k+1}\|}{\|E_k\|}\right)^{(2)} = \frac{\left\|E_{k+1}^{(1)} + \frac{\Delta E_{k+1}^{(1)}}{1-\lambda_1}\right\|}{\left\|E_k^{(1)} + \frac{\Delta E_k^{(1)}}{1-\lambda_1}\right\|} \quad\quad\quad \ldots\ldots\ldots\ldots [13]$$

In terms of the Eigen values $\lambda$ and Eigen vectors $v$ of the iteration matrix M, the terms in the norm expression of equation (13) above are given by:

$$E_k^{(1)} = c_1\lambda_1^{k}v_1 + c_2\lambda_2^{k}v_2 + c_3\lambda_3^{k}v_3 + \ldots\ldots + c_N\lambda_N^{k}v_N$$

$$E_{k+1}^{(1)} = c_1 \lambda_1^{k+1} v_1 + c_2 \lambda_2^{k+1} v_2 + c_3 \lambda_3^{k+1} v_3 + \ldots\ldots + c_N \lambda_N^{k+1} v_N$$

$$\Delta E_k^{(1)} = E_{k+1}^{(1)} - E_k^{(1)}$$

$$= c_1(\lambda_1 - 1)\lambda_1^k v_1 + c_2(\lambda_2 - 1)\lambda_2^k v_2 + c_3(\lambda_3 - 1)\lambda_3^k v_3 + \ldots\ldots + c_N(\lambda_n - 1)\lambda_N^k v_N$$

$$\Delta E_{k+1}^{(1)} = E_{k+2}^{(1)} - E_{k+1}^{(1)}$$

$$= c_1(\lambda_1 - 1)\lambda_1^{k+1} v_1 + c_2(\lambda_2 - 1)\lambda_2^{k+1} v_2 + c_3(\lambda_3 - 1)\lambda_3^{k+1} v_3 + \ldots\ldots + c_N(\lambda_n - 1)\lambda_N^{k+1} v_N$$

Collecting the terms for both the numerator and denominator yields,

$$\frac{\left\|E_{k+1}^{(1)} + \frac{\Delta E_{k+1}^{(1)}}{1 - \lambda_1}\right\|}{\left\|E_k^{(1)} + \frac{\Delta E_k^{(1)}}{1 - \lambda_1}\right\|} = \frac{\left\|\sum_{i=1}^{N} c_i (\lambda_i)^{k+1} v_i \left[1 - \frac{\lambda_i - 1}{\lambda_1 - 1}\right]\right\|}{\left\|\sum_{i=1}^{N} c_i (\lambda_i)^k v_i \left[1 - \frac{\lambda_i - 1}{\lambda_1 - 1}\right]\right\|} \quad \ldots\ldots\ldots\ldots [14]$$

Examination of the terms on the right hand side of equation (14) above reveals that the first terms of the summation (i.e. i=1) in both the numerator and denominator vanish. Therefore, the second-order Aitken extrapolation reduces the Eigenvalue so that $\lambda_2$ becomes the dominant Eigen value. As will be shown below, the second dominant Eigen value of the iteration matrix, $\lambda_2$, will be equal to the ratio of the error vectors for the second order Aitken extrapolation.

Factoring out the $\lambda_2$ term in equation (14) above will result in the following expression:

$$\frac{\left\|E_{k+1}^{(1)} + \frac{\Delta E_{k+1}^{(1)}}{1 - \lambda_1}\right\|}{\left\|E_k^{(1)} + \frac{\Delta E_k^{(1)}}{1 - \lambda_1}\right\|} = \frac{\lambda_2^{k+1} \left\|c_2 v_2 \left[1 - \frac{\lambda_2 - 1}{\lambda_1 - 1}\right] + \sum_{i=3}^{N} c_i \left(\frac{\lambda_i}{\lambda_2}\right)^{k+1} v_i \left[1 - \frac{\lambda_i - 1}{\lambda_1 - 1}\right]\right\|}{\lambda_2^k \left\|c_2 v_2 \left[1 - \frac{\lambda_2 - 1}{\lambda_1 - 1}\right] + \sum_{i=3}^{N} c_i \left(\frac{\lambda_i}{\lambda_2}\right)^k v_i \left[1 - \frac{\lambda_i - 1}{\lambda_1 - 1}\right]\right\|}$$

Once $\lambda_1$ has been eliminated and $\lambda_2$ is the dominant Eigen value, the ratios:

$$\left(\frac{\lambda_i}{\lambda_2}\right)^k, \left(\frac{\lambda_i}{\lambda_2}\right)^{k+1} \quad i = 3, 4, 5 \ldots\ldots$$

being less than one, will vanish at higher k values so that the error norm ratios become;

$$\left(\frac{\|E_{k+1}\|}{\|E_k\|}\right)^{(2)} = \frac{\left\|E_{k+1}^{(1)} + \frac{\Delta E_{k+1}^{(1)}}{1-\lambda_1}\right\|}{\left\|E_k^{(1)} + \frac{\Delta E_k^{(1)}}{1-\lambda_1}\right\|} = \frac{\lambda_2^{k+1}\left\|c_2 v_2\left[1 - \frac{\lambda_2 - 1}{\lambda_1 - 1}\right]\right\|}{\lambda_2^k\left\|c_2 v_2\left[1 - \frac{\lambda_2 - 1}{\lambda_1 - 1}\right]\right\|}$$

$$\left(\frac{\|E_{k+1}\|}{\|E_k\|}\right)^{(2)} = \lambda_2$$

Similarly, it is easy to show that for the third and higher order Aitken extrapolations the error vector ratios correspond to the i[th] Eigen value of the iteration matrix. The higher order decomposition of Eigen values through higher order Aitken extrapolation can be generalized as:

$$\left(\frac{\|E_{k+1}\|}{\|E_k\|}\right)^{(i)} = \lambda_i \quad for\ i = 1, 2, 3, \ldots \ldots \ldots N \qquad \text{............. [15]}$$

In effect, higher order Aitken extrapolation successively decomposes the dominant Eigen values so that the error terms are determined eventually by the lowest Eigen value of the iteration matrix. This works for both the convergent and divergent Gauss-Seidel iteration. However, the procedure is best in decomposing the first two dominant Eigen values beyond which the decomposition might be slow or inexact due to the successively small difference in the error vectors. The example that follow later for diverging Gauss-Seidel iteration illustrate this fact.

**Coupling of SOR technique with Geometric series extrapolation**

The extrapolation to the Gauss-Seidel iteration can well be extended to the successive over relaxation (SOR) method. In matrix form, the SOR iteration process is:

$$X^{(k+1)} = -(D + \omega L)^{-1} * [\omega U + (\omega - 1)D]X^{(k)} + (D + \omega L)^{-1} * \omega B \qquad \text{............. [16]}$$

Where ω is the relaxation factor and the other terms are as deifned above. The iteration matrix is the coefficient of the $X^k$ term in equation (16) above and is given by:

$$M_\omega = -(D + \omega L)^{-1} * [\omega U + (\omega - 1)D] \quad\quad \text{............ [17]}$$

The iteration formula for the successive over relaxation technique in terms of individual x values is given by;

$$x_i^{k+1} = (1 - \omega)x_i^{(k)} + \frac{\omega}{a_{ii}}\left(b_i - \sum_{j<i} a_{ij}x_j^{(k+1)} - \sum_{j>i} a_{ij}x_j^{(k)}\right) \quad\quad \text{............ [18]}$$

$$i = 1,2,\ldots.n$$

The acceleration factor ω cannot be easily determined in advance. It depends on the coefficient matrix A. If the coefficient matrix A is symmetric as well as positive definite, the spectral radius of the iteration matrix $M_\omega$ will be less than one - ensuring convergence of the process - when the ω value lies between 0 and 2.

The procedure for extrapolation of the SOR process using geometric series sum, based on the dominant Eigen value of the iteration matrix as a ratio of the geometric series, follows a similar process to the one mentioned earlier. The only change is in the iteration matrix which is modified by the relaxation factor ω while the condition for convergence (i.e. the dominant Eigen value of the iteration matrix being less than one) remains the same.

However, it should be noted that the optimum relaxation factor ω is not necessarily the same as the SOR – optimum when the SOR technique is combined with Aitken extrapolation. This is illustrated in the application example of the heat flow problem presented in this paper. The $\omega_{opt}$ for the SOR techniques is so chosen that the two dominant Eigen-values become equal in magnitude and these Eigen values at the optimum ω can be complex numbers leading possibly to the failure of extrapolation methods. The fact that, at the optimum value of the acceleration factor ω, the dominant Eigen values turn out to be complex numbers is also shown in the heat flow example presented in this paper. Coupling the SOR technique with the Aitken extrapolation at the exact optimum ω is not necessary as the result is not very sensitive to the ω value as will be shown in the heat flow example that follows. However, failure is not necessarily always the case for coupling at the optimum ω value. The application example suggests that in the case of coupling of SOR with Aitken's extrapolation, the ω value can be chosen so that it is slightly less than the SOR optimum value enabling the coupling to be made at real Eigen values without the method failing to lead to convergence to the solution. .

## 4. Examples from convergent iterations

*Example 1:*

The first example below is a simple 2 x 2 equation in x and y values.

$$2x + y = 7$$
$$x - y = 2$$

The Gauss – Seidel iteration starts at values distant from the solution, i.e. the starting values of x = 10,000 and y = 4250. The x and y values of the iteration, differences in consecutive steps of the iteration $e_i$ and the ratios, $\lambda$, are computed and shown in Table 1 below

Table 1 Computation for geometric series extrapolation of the Gauss-Seidel iteration for the 2 by 2 equation

| Iteration step | x | y | Consecutive Difference $e_i$ (x) | Consecutive Difference $e_i$ (y) | Ratio $e_{i+1}/e_i$ (x) | Ratio $e_{i+1}/e_i$ (y) |
|---|---|---|---|---|---|---|
| 0 | 10000.00 | 4250.00 | -12121.50 | -6373.50 | -0.26 | -0.50 |
| 1 | -2121.50 | -2123.50 | 3186.75 | 3186.75 | -0.50 | -0.50 |
| 2 | 1065.25 | 1063.25 | -1593.38 | -1593.38 | -0.50 | -0.50 |
| 3 | -528.13 | -530.13 | 796.69 | 796.69 | -0.50 | -0.50 |
| 4 | 268.56 | 266.56 | -398.34 | -398.34 | -0.50 | -0.50 |
| 5 | -129.78 | -131.78 | 199.17 | 199.17 | -0.50 | -0.50 |
| 6 | 69.39 | 67.39 | -99.59 | -99.59 | -0.50 | -0.50 |
| 7 | -30.20 | -32.20 | 49.79 | 49.79 | -0.50 | -0.50 |
| 8 | 19.60 | 17.60 | -24.90 | -24.90 | -0.50 | -0.50 |
| 9 | -5.30 | -7.30 | 12.45 | 12.45 | -0.50 | -0.50 |
| 10 | 7.15 | 5.15 | -6.22 | -6.22 | -0.50 | -0.50 |
| 11 | 0.93 | -1.07 | 3.11 | 3.11 | -0.50 | -0.50 |
| 12 | 4.04 | 2.04 | -1.56 | -1.56 | -0.50 | -0.50 |
| 13 | 2.48 | 0.48 | 0.78 | 0.78 | -0.50 | -0.50 |
| 14 | 3.26 | 1.26 | -0.39 | -0.39 | | |
| 15 | 2.87 | 0.87 | -2.87 | -0.87 | | |

It is clear that the ratios of differences in x and y values converge almost immediately to -0.50 for both the x and y values of the iteration. It will be later shown that, this value is the largest Eigen-value of the iteration matrix, M.

For calculating the extrapolated x value of the iteration, x= 10,000 cannot be used as the $x_0$ value since the ratio at this level (-0.26) is not sufficiently convergent (compared to -0.5). Therefore, the second x value, i.e., $x_1$ = -2125.5 is taken. For this value the $e_{x0}$ value is listed in Table 1 as 3186.75. For the y iteration $y_0$ = 4250 can be taken as the ratio converged immediately with the first iteration. The $e_{yo}$ value is also taken to be -6373.5

The $x_\infty$ and $y_\infty$ values are calculated as;

$$x_\infty = x_0 + \frac{e_{x0}}{1-\lambda_x} = -2121.5 + \frac{3186.75}{1-(-0.5)} = 3.000$$

$$y_\infty = y_0 + \frac{e_{y0}}{1-\lambda_y} = 4250 + \frac{-6373.5}{1-(-0.5)} = 1.000$$

These values as expected are exactly equal to the solution of the equations.

To see that the ratios of alternative differences are the same as the largest Eigen value of the iteration matrix, the equation;

$$2x + y = 7$$
$$x - y = 2$$

is written as AX = B;

$$\begin{bmatrix} 2 & 1 \\ 1 & -1 \end{bmatrix} \begin{bmatrix} x \\ y \end{bmatrix} = \begin{bmatrix} 7 \\ 2 \end{bmatrix}$$

Using the relation A = L + D + U

$$A = \begin{bmatrix} 2 & 1 \\ 1 & -1 \end{bmatrix} \; ; \; L = \begin{bmatrix} 0 & 0 \\ 1 & 0 \end{bmatrix} \; ; \; D = \begin{bmatrix} 2 & 0 \\ 0 & -1 \end{bmatrix} \; ; \; U = \begin{bmatrix} 0 & 1 \\ 0 & 0 \end{bmatrix}$$

It can be shown that using the L, D and U matrices;

$$(L+D)^{-1}U = \begin{bmatrix} 0 & 1/2 \\ 0 & 1/2 \end{bmatrix}$$

and that;

$$M = -[(L+D)^{-1}U] = \begin{bmatrix} 0 & -1/2 \\ 0 & -1/2 \end{bmatrix}$$

The Eigen values of the iteration matrix M matrix are found by solving the determinant equation;

$$\det(M - \lambda I) = \begin{vmatrix} 0 - \lambda & -1/2 \\ 0 & -\frac{1}{2} - \lambda \end{vmatrix} = 0$$

$$(-\lambda) * \left(-\frac{1}{2} - \lambda\right) = 0$$

$$\lambda = 0 \text{ or } \lambda = -\frac{1}{2} = -0.5$$

The largest Eigen value is -0.5 and it can be seen from Table 1 that the ratios of consecutive differences for both x and y variables converge to the largest Eigen value of the iteration matrix.

Example 2

The 2$^{nd}$ example is a 3 x 3 systems of linear equations given below

$$3x + y + z = 5$$
$$-x + 2 + 3z = 4$$
$$x + y - 4z = -2$$

The solution vector is $[x, y, z] = [1, 1, 1]^T$

The iteration starts with the vector : $[x, y, z]^T = [10, -8, 5]^T$

The Gauss-Seidel iteration was carried out 13 times at which level four decimal-digit accuracy was obtained for the ratio of consecutive differences in x, y and z values. Table 2 below shows the corresponding x, e values and ratios of consecutive e values.

Table 2 Results geometric series extrapolation of the Gauss-Seidel iteration for the 3 by 3 equation after 9 iterations.

| Values after 9 iterations | Differences in values | Ratios of differences | Extrapolated values |
|---|---|---|---|
| $x_0$ = 1.44846653 | $e_{x0}$ = -0.81586443 | $\lambda_x$ = -0.81923923 | $x_\infty$ = 1.000001908 |
| $y_0$ = 1.79206166 | $e_{y0}$ = -1.44095868 | $\lambda_y$ = -0.81924816 | $y_\infty$ = 0.999998918 |
| $z_0$ = 1.31013205 | $e_{z0}$ = -0.56420578 | $\lambda_z$ = -0.81924493 | $z_\infty$ = 1.000000209 |

It is shown in Table 2 that, with the geometric series extrapolation a 5 digit accuracy has been obtained for the solution with just 9 iterations. The normal Gauss-Seidel iteration requires more than 60 iterations to arrive at 5 digit accuracy.

For the coefficient matrix A of the given equation, it can be shown that the iteration matrix M is given by;

$$M = -[(L+D)^{-1}U] = \begin{bmatrix} 0 & -\frac{1}{3} & -\frac{1}{3} \\ 0 & -\frac{1}{6} & -\frac{10}{6} \\ 0 & -\frac{1}{8} & -\frac{1}{2} \end{bmatrix}$$

The Eigen values are found by solving the determinant;

$$\det(M - \lambda I) = \begin{vmatrix} 0-\lambda & -\frac{1}{3} & -\frac{1}{3} \\ 0 & -\frac{1}{6}-\lambda & -\frac{10}{6} \\ 0 & -\frac{1}{8} & -\frac{1}{2}-\lambda \end{vmatrix} = 0$$

$$(-\lambda) * \left[\lambda^2 + \frac{2}{3}\lambda - \frac{1}{8}\right] = 0$$

$$\lambda = 0, \quad \lambda = 0.152579324 \quad or \quad \lambda = -0.81924599$$

Therefore, the largest Eigen value of the iteration matrix M is $\lambda = -0.81924599$ and all the ratios used above were approaching towards the largest Eigen value of the iteration matrix accurate to 5 digits as shown in Table 2.

*Application Example (Heat flow problem)*

Example of application of the Aitken extrapolation method for a steady state heat flow problem involving Laplace equation is given below [15]. A rectangular thin steel plate with dimension 10 x 20 cm has one of its 10 cm edges held at 100 $^0c$ and the other three edges are held at $0^0c$. The thermal conductivity is given as k = 0.16 cal/sec.cm$^2$.$^0$c/cm. Figure 1 below shows the steel plate steady state conditions temperatures.

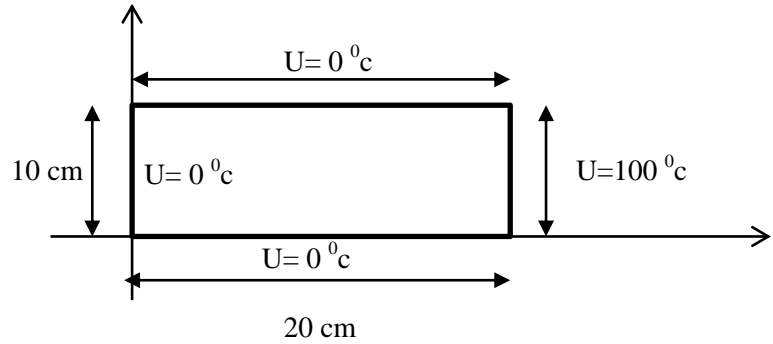

Figure 1  Rectangular plate for the steady state heat flow problem with boundary conditions.

The steady state heat flow problem is described by the Laplace equation:

$$\frac{\partial^2 u}{\partial x^2} + \frac{\partial^2 u}{\partial y^2} = 0$$

With the boundary conditions, $u(x,0) = u(x,10) = u(0,y) = 0$ and $u(20,y) = 100^0C$.

The nine-point finite difference formula of the Laplace equation is used for computation. This formula is symbolically represented by:

$$\nabla^2 u_{i,j} = \frac{1}{6h^2} \begin{Bmatrix} 1 & 4 & 1 \\ 4 & -20 & 4 \\ 1 & 4 & 1 \end{Bmatrix} u_{i,j} = 0$$

The algebraic form of the difference equation is:

$$\frac{1}{6h^2} [u_{i+1,j+1} + u_{i+1,j-1} + u_{i-1,j+1} + u_{i-1,j-1} + 4u_{i+1,j} + 4u_{i-1,j} + 4u_{i,j+1} + 4u_{i,j-1} - 20\, u_{i,j}] = 0$$

Where h is the grid size and the u values are the temperatures at the intersection of the rectangular grid lines.

Using a grid size of 2.5 cm, the 21 interior grid points shown in Figure 2 below are generated.

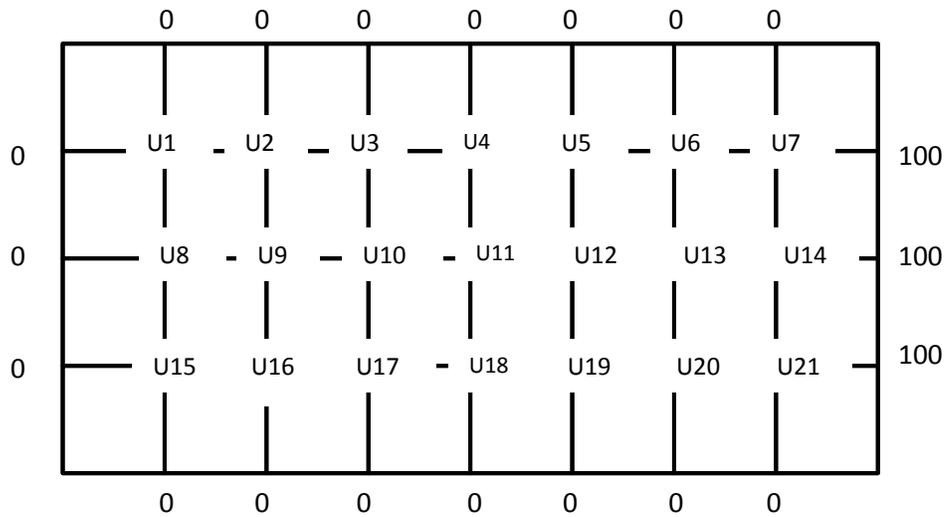

Figure 2   Rectangular plate for the heat flow problem with interior points and boundary conditions shown

The coefficient matrix A of the linearized form of the Laplace equation AU= B by finite differencing is given in Table 3 below, followed by the right hand side vector B.

A MATLAB program was written to solve the linear system of equations for the 21 unknown temperatures using the Gauss-Seidel iteration initially. In addition the geometric series extrapolation of the Gauss-Seidel was computed by writing the corresponding program in the MATLAB environment.  The solution vector X for each of the 100 iterations for the normal Gauss-Seidel iteration was computed. In addition the Euclidian norms of the error vectors for each step of the iteration were computed.  For the geometric series extrapolation, in addition to the solution vectors and norm of the error vectors, the consecutive difference ratios, as approximation of the maximum Eigen-value were also computed.

Fig 3 below shows a comparison of the number of iterations required to arrive at more or less the same magnitude of the norm of the error vector for the normal Gauss-Seidel iteration and for the extrapolation. It is observed from Figure 3 that the extrapolation procedure will give an acceleration factor in the range between 2.1 and 2.2. For example whereas 16 iterations are required for the normal Gauss-Seidel iteration  giving error norm of $0.06^+$ , the geometric series extrapolation required only 7 iterations. For the subsequently smaller error norms, the numbers of iterations required are in the ratio of   16/7, 18/8, 20/9, 22/10, 24/11.  It is clear that as the error norm gets smaller, the acceleration factor approaches a value of 2.

**Table 3 coefficient matrix A of the finite difference form of the heat flow problem**

```
-20  4   0   0   0   0   0   4   1   0   0   0   0   0   0   0   0   0   0   0   0
 4  -20  4   0   0   0   0   1   4   1   0   0   0   0   0   0   0   0   0   0   0
 0   4  -20  4   0   0   0   0   1   4   1   0   0   0   0   0   0   0   0   0   0
 0   0   4  -20  4   0   0   0   0   1   4   1   0   0   0   0   0   0   0   0   0
 0   0   0   4  -20  4   0   0   0   0   1   4   1   0   0   0   0   0   0   0   0
 0   0   0   0   4  -20  4   0   0   0   0   1   4   1   0   0   0   0   0   0   0
 0   0   0   0   0   4  -20  0   0   0   0   0   1   4   0   0   0   0   0   0   0
 4   1   0   0   0   0   0  -20  4   0   0   0   0   0   4   1   0   0   0   0   0
 1   4   1   0   0   0   0   4  -20  4   0   0   0   0   1   4   1   0   0   0   0
 0   1   4   1   0   0   0   0   4  -20  4   0   0   0   0   1   4   1   0   0   0
 0   0   1   4   1   0   0   0   0   4  -20  4   0   0   0   0   1   4   1   0   0
 0   0   0   1   4   1   0   0   0   0   4  -20  4   0   0   0   0   1   4   1   0
 0   0   0   0   1   4   1   0   0   0   0   4  -20  4   0   0   0   0   1   4   1
 0   0   0   0   0   1   4   0   0   0   0   0   4  -20  0   0   0   0   0   1   4
 0   0   0   0   0   0   0   4   1   0   0   0   0   0  -20  4   0   0   0   0   0
 0   0   0   0   0   0   0   1   4   1   0   0   0   0   4  -20  4   0   0   0   0
 0   0   0   0   0   0   0   0   1   4   1   0   0   0   0   4  -20  4   0   0   0
 0   0   0   0   0   0   0   0   0   1   4   1   0   0   0   0   4  -20  4   0   0
 0   0   0   0   0   0   0   0   0   0   1   4   1   0   0   0   0   4  -20  4   0
 0   0   0   0   0   0   0   0   0   0   0   1   4   1   0   0   0   0   4  -20  4
 0   0   0   0   0   0   0   0   0   0   0   0   1   4   0   0   0   0   0   4  -20
```

The right hand side vector of the finite difference equation is given as:

$B = [0\ 0\ 0\ 0\ 0\ 0\ -550\ 0\ 0\ 0\ 0\ 0\ 0\ -600\ 0\ 0\ 0\ 0\ 0\ 0\ -550\ ]^T$

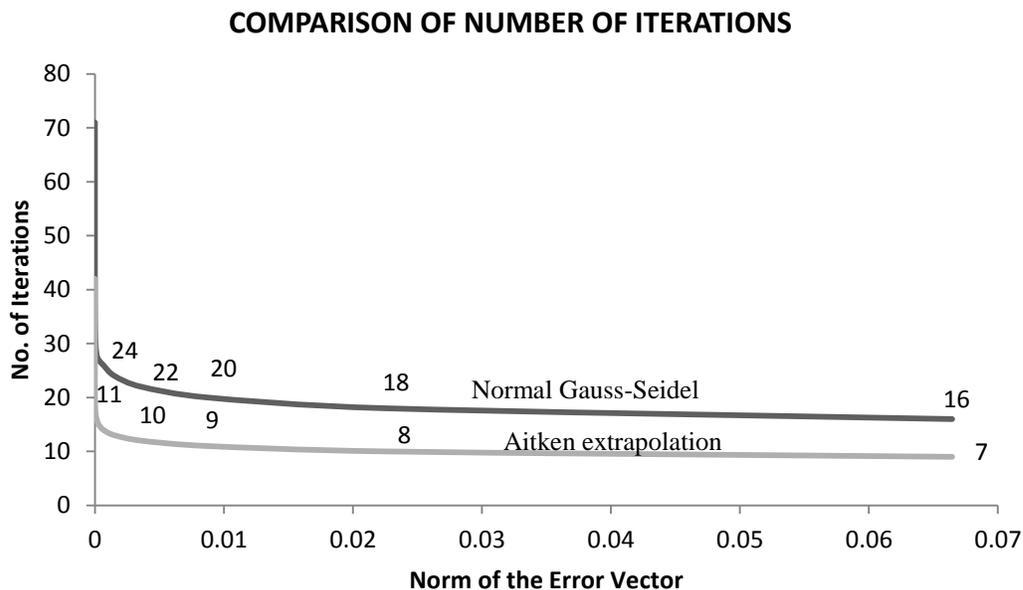

Figure 3   Comparison of number of iterations required for given norm of error vector between the normal Gauss-Seidel iteration and the geometric series extrapolation

*Application Example - Double acceleration combined with the SOR technique*

For the rectangular plate heat flow problem given above with 32 mesh divisions of interval Δh = 2.5 cm in both x and y directions forming 21 interior points for the Gauss-Seidel iteration, the successive over relaxation technique (SOR) was applied in conjunction with the geometric series sum based Aitken's extrapolation. The optimum acceleration factor ω to be used in equations (17) and (18) for rectangular regions with uniform boundary conditions (Drichlet conditions) is given by [15];

$$\omega_{opt} = \frac{4}{2 + \sqrt{4 - c^2}} \qquad \text{............ [19]}$$

The value of c in equation (19) above is given by;

$$c = \left(\cos\frac{\pi}{p} + \cos\frac{\pi}{q}\right)$$

Where p and q are the number of mesh divisions in the x and y directions. For the given problem p = 8 and q = 4, so that;

$$c = \left(\cos\frac{\pi}{8} + \cos\frac{\pi}{4}\right) = 1.63098$$

$$\omega_{opt} = \frac{4}{2 + \sqrt{4 - 1.63098^2}} = 1.267$$

Therefore, the minimum number of iterations required is achieved by using the optimum acceleration factor ω of 1.267. This is shown in Table 5 below for the SOR column where the minimum number of iteration of 35 was required to reach to the solution vector within $10^{-15}$ accuracy. The normal Gauss-Seidel process required 80 iterations whereas the geometric series extrapolation needed 47 iterations. Coupling of the SOR technique with geometric series extrapolation resulted in further reduction in the number of iterations required from 35 to 24, a reduction by about 31 % from the SOR result.

Table 5. Comparison of number of iterations required between SOR and Aitken's geometric series extrapolation coupled with the SOR technique.

| w | SOR+GSE* | SOR | Extra acceleration | Normal | Largest Eigen value of the SOR iteration matrix |
|---|---|---|---|---|---|
| 0.8 | 78 | 127 | 1.63 | 80 | 0.859905 |
| 0.9 | 59 | 102 | 1.73 | 80 | 0.827175 |
| 1 | 47 | 80 | 1.70 | 80 | 0.788581 |
| 1.05 | 40 | 70 | 1.75 | 80 | 0.760252 |
| 1.1 | 34 | 62 | 1.82 | 80 | 0.729593 |
| 1.15 | 26 | 54 | 2.08 | 80 | 0.691062 |
| 1.2 | 25 | 45 | 1.80 | 80 | 0.637938 |
| 1.23 | 24 | 40 | 1.67 | 80 | 0.59002 |
| 1.25 | 26 | 35 | 1.35 | 80 | 0.532422 |
| 1.267 | 27 | 35 | 1.30 | 80 | Complex number |
| 1.24 | 25 | 37 | 1.48 | 80 | Complex number |
| 1.3 | 32 | 37 | 1.16 | 80 | Complex number |
| 1.4 | 41 | 48 | 1.17 | 80 | Complex number |
| 1.6 | 72 | 76 | 1.06 | 80 | Complex number |
| 1.8 | 298 | 168 | 0.56 | 80 | Complex number |

* Coupling of SOR technique with Aitken extrapolation

The rate of reduction in number of iterations required by coupling SOR with Aitken geometric series extrapolation is displayed in Figure 4 below for different values of the SOR acceleration factor, ω based on the values given in Table 5 above. It can be seen from the figure that significant reduction in the number of iterations required is achieved for ω values between 1 and $\omega_{opt}$. In fact, the optimum value of ω for the coupled iteration (SOR + Aitken geometric series extrapolation) lies below the SOR optimum for ω, i.e., at ω = 1.23.

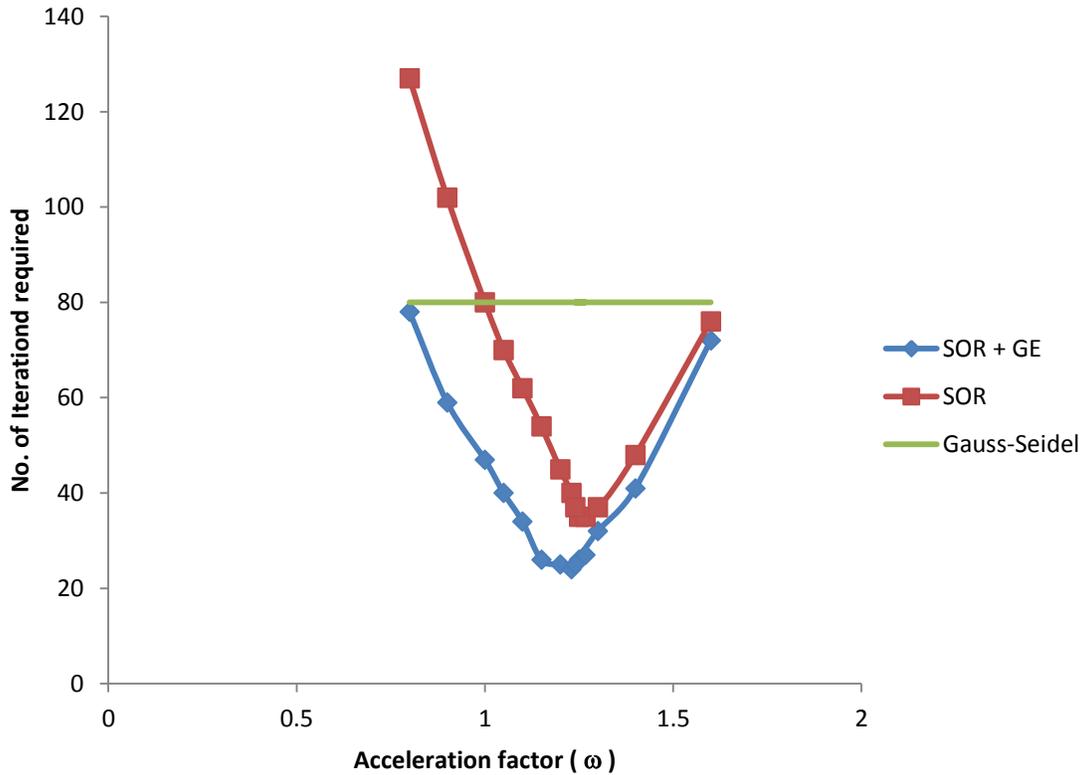

Figure 4  Comparison of number of iterations required for given norm of error vector between the SOR method and the geometric series extrapolation coupled with SOR.
In the figure, GE = Geometric series extrapolation. SOR = Successive Over Relaxation.

Beyond the optimum acceleration factor, the difference between SOR and the coupled iteration is insignificant suggesting that a value ω slightly less than the SOR optimum should be used if coupling is to be made. It is interesting to note that the largest Eigen values of the iteration matrix for the SOR technique, i.e.,

$$M_\omega = (-(D + \omega L)^{-1} * [\omega U + (\omega - 1)D])$$

turn out to be a complex numbers at the optimum ω value and beyond (refer to table 5 above for ω = 1.267 and above). For all the ω values above the SOR optimum, the corresponding largest Eigen values are complex numbers and the spectra radii are increasing. The progressive reduction in largest Eigen values of the iteration matrix is also evident as the ω value increases towards the SOR optimum. However, the extrapolation did not fail even if the dominant Eigen values were complex numbers at end beyond the optimum ω values.

The advantage of coupling the SOR technique with Aitken extrapolation is evident from the above example. In addition, the Aitken extrapolation is based on the Gauss-Seidel iteration and as such does not involve extra calculation except generating a geometric series sum. The optimum SOR value is not always predictable. However, Aitken's geometric series extrapolation can still work with or without the use of optimum ω value. The above example shows that for acceleration factors slightly greater than one, significant reduction in the number of iterations required were obtained when the SOR was coupled with extrapolation.

## 5. Examples from a Divergent Gauss-Seidel Iteration

*Example 5.1 A 2 X 2 equation with diagonally non-dominant coefficient matrix.*

The first example below is a simple 2 x 2 system of equations with diagonally non-dominant coefficient matrix

$$2x + 10y = 12$$
$$15x - 5y = 10$$

The Gauss – Seidel iteration starts at values of x = 8 and y = 10. The x and y values of the iteration, differences in consecutive steps of the iteration $e_i$ and the ratios, $\lambda$, were computed and are shown in Table 6 below. As Table 6 shows, the Gauss-Seidel based iteration diverges as the coefficient matrix is also diagonally non-dominant. This is also evident from the ratio of differences in x and y values shown in Table 6. This ratio is -15 for both x and y iterations and it to corresponds to the dominant Eigen value of the iteration matrix. However, the Aitken extrapolation iterations (shown in the last columns of the table) invariably converge to the true solutions x=1 and y=1.

Using the values of the second iteration in Table 6, The extrapolated x and y values, $x_s$ and $y_s$ are calculated using equation (12) as;

$$x_s = x_0 + \frac{e_{x0}}{1 - \lambda_x} = 676 + \frac{-10800}{1 - (-15)} = 1.000$$

$$y_s = y_0 + \frac{e_{y0}}{1 - \lambda_y} = 2026 + \frac{-32400}{1 - (-15)} = 1.000$$

These values as expected are exactly equal to the solution of the equations.

Table 6. Gauss-Seidel iteration results for a diagonally non-dominant and diverging 2x2 system of equations.

| Iteration step | x | y | Consecutive Difference $e_i$ (x) | Consecutive Difference $e_i$ (y) | Ratio $e_{i+1}/e_i$ (x) | Ratio $e_{i+1}/e_i$ (y) | Extrapolation (X) | Extrapolation (Y) |
|---|---|---|---|---|---|---|---|---|
| 0 | 8 | 10 | -52 | -144 | -13.846 | -15 | 4.49 | 1 |
| 1 | -44 | -134 | 720 | 2160 | -15 | -15 | 1 | 1 |
| 2 | 676 | 2026 | -10800 | -32400 | -15 | -15 | 1 | 1 |
| 3 | -10124 | -30374 | 162000 | 486000 | -15 | -15 | 1 | 1 |
| 4 | 151876 | 455626 | -2430000 | -7290000 | -15 | -15 | 1 | 1 |
| 5 | -2278124 | -6834374 | 36450000 | $1.09 \times 10^{08}$ | -15 | -15 | 1 | 1 |
| 6 | 34171876 | 102515626 | $-5.5E+^{08}$ | $-1.6 \times 10^{09}$ | -15 | -15 | 1 | 1 |
| 7 | $-5.1 \times 10^{08}$ | $-1.538 \times 10^{09}$ | $8.2 \times 10^{09}$ | $2.46 \times 10^{10}$ | -15 | -15 | 1 | 1 |
| 8 | $7.69 \times 10^{*09}$ | $2.3066 \times 10^{10}$ | $-1.2 \times 10^{11}$ | $-3.7 \times 10^{11}$ | -15 | -15 | 1 | 1 |
| 9 | $-1.2 \times 10^{11}$ | $-3.46 \times 10^{11}$ | $1.85 \times 10^{12}$ | $5.54 \times 10^{12}$ | | | | |
| 10 | $1.73 \times 10^{12}$ | $5.1899 \times 10^{12}$ | $-1.7 \times 10^{12}$ | $-5.2 \times 10^{12}$ | | | | |

*Example 5.2 A 4X 4 equations with diagonally non-dominant coefficient matrix.*

A second example of a four by four system of equations in which the coefficient matrix is once again diagonally non-dominant is presented below.

$$\begin{bmatrix} 20 & 234 & 123 & 20 \\ 136 & 56 & 120 & 125 \\ 123 & 120 & 76 & 25 \\ 20 & 125 & 145 & 20 \end{bmatrix} \begin{bmatrix} x_1 \\ x_2 \\ x_3 \\ x_4 \end{bmatrix} = \begin{bmatrix} -783 \\ -346 \\ -127 \\ -481 \end{bmatrix}$$

The normal Gauss-Seidel iteration quickly diverges as expected and it is not possible to reach to the solution as such. Higher order Aitken extrapolation has to be carried out since the two dominant Eigen values of the iteration matrix are large (16.7 and 5.77 respectively). A fifth order Aitken extrapolation successfully achieves convergence whereas the first two dominant Eigen values were successfully decomposed with the first and second order Aitken extrapolations respectively.

The solution of the system of equations together with the Eigen values of the iteration matrix for the Gauss-Seidel iteration as obtained using MATLAB program are given in Table 7 below.

Table 7. Values of solution to the four by four system of equations and Eigen values of the iteration matrix using MATLAB program.

| Solutions of the systems of equation | Eigen values of the iteration matrix |
|---|---|
| $X_1 = 3.054225004761563$ | -16.700300071436452 |
| $X_2 = -2.904223059942874$ | 5.774221605390342 |
| $X_3 = -0.661832433353327$ | 0.085523954767930 |
| $X_4 = -4.154545738306979$ | 0 |

Figure 5 below shows the decomposition of error ratios for the five-order Aitken extrapolation. Comparing the ratio of successive errors given at each level of extrapolation with the Eigen values of the iteration matrix given in Table 7 above shows that the first two dominant Eigen values are decomposed exactly ( i.e. $\lambda_1 = -16.7003$ and $\lambda_2 = 5.7742$). However, for the $3^{rd}$ and higher order extrapolation the decomposition slows but eventually reduces sufficiently enabling convergence of the iteration.

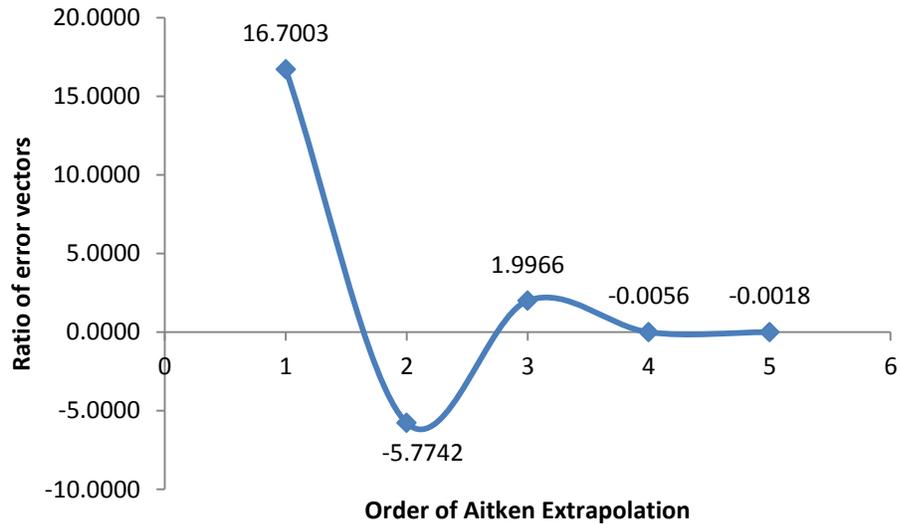

Figure 5. Variation of ratios of the error vectors at $1^{st}$ to $5^{th}$ order Aitken extrapolation

Figure 6 below shows the progress of reduction of the norm of the error vector computed from

$$\|E\| = \|B - AX\|$$

As can be seen from Table 8, with the fifth order Aitken extrapolation, the norm of the error vector eventually reduces down to $10^{-9}$. It should be recalled that the normal Gauss-Seidel iteration is rapidly diverging for this system of equations. On the other hand higher order Aitken extrapolation as applied in this example successfully converges to the solution of the system of equations.

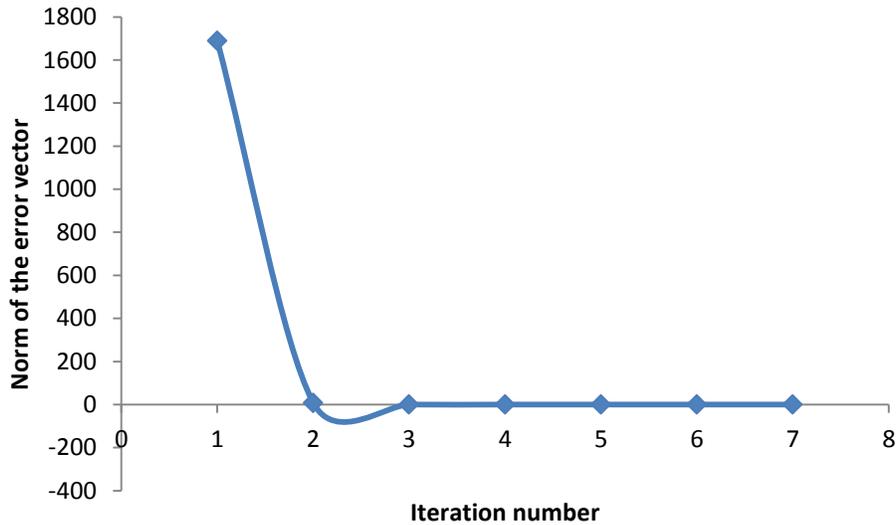

Figure 6. Progress of reduction of the error vector for a $5^{th}$ order Aitken extrapolation

Table 8. Progress of the $5^{th}$ order Aitken extrapolation for the 4 x 4 system of equations

| Iteration number | $X_1$ | $X_2$ | $X_3$ | $X_4$ | Norm of the error vector |
|---|---|---|---|---|---|
| 1 | 1 | 1 | 1 | 1 | 1689.086 |
| 2 | 3.044218606 | -2.904617307 | -0.622362121 | -4.153534798 | 8.14532 |
| 3 | 3.05422335 | -2.904224136 | -0.66182978 | -4.15448631 | 0.00781 |
| 4 | 3.054225007 | -2.904223059 | -0.66183244 | -4.154547438 | 0.000223 |
| 5 | 3.054225005 | -2.90422306 | -0.661832433 | -4.15454574 | $6.37*10^{-06}$ |
| 6 | 3.054225005 | -2.90422306 | -0.661832433 | -4.154545738 | $1.83*10^{-07}$ |
| 7 | 3.054225005 | 2.90422306 | -0.661832433 | -4.154545738 | $5.32*10^{-09}$ |

*Example 5.3 A 6X 6 system of equations with a diagonally non-dominant coefficient matrix.*

A further example of diagonally non-dominant six by six systems of linear equations is given below:

$$\begin{bmatrix} 400 & 35 & 432 & 10 & 4820 & 0 \\ 35 & 600 & 485 & 30 & 20 & 2000 \\ 196 & 10 & 545 & 48 & 34 & 974 \\ 0 & 30 & 48 & 631 & 20 & 347 \\ 4820 & 20 & 34 & 545 & 768 & 0 \\ 0 & 2000 & 800 & 0 & 874 & 657 \end{bmatrix} \begin{bmatrix} x1 \\ x2 \\ x3 \\ x4 \\ x5 \\ x6 \end{bmatrix} = \begin{bmatrix} -7830 \\ -3460 \\ -1270 \\ -4810 \\ -7623.8 \\ -2690 \end{bmatrix}$$

The dominant Eigen value of the Gauss-Seidel iteration matrix is -75.7966. Table 9 below shows the exact solutions of the system of equations together with the Eigen values of the iteration matrix which were obtained from a MATLAB program. With the combination of dominant Eigen values shown in Table 9, the normal Gauss-Seidel iteration quickly diverges. However, a $4^{th}$ order Aitken extrapolation was enough to bring the iteration to convergence.

Table 9. Values of solution to the six by six system of equations and Eigen values of the iteration matrix.

| Solutions of the systems of equation | Eigen values of the iteration matrix |
| --- | --- |
| $X_1 = $ -0.563147393304287 | -75.796638515975403 |
| $X_2 = $ -0.731832157439239 | -11.704130560629785 |
| $X_3 = $ 1.857839885254053 | -0.368124943597328 |
| $X_4 = $ -6.666186315796603 | -0.027943199473836 |
| $X_5 = $ -1.725114498846410 | -0.000000000000001 |
| $X_6 = $ -1.833877505834098 | 0 |

Figure 7 below shows the decomposition of the ratio of error vectors at each order of Aitken extrapolation. At the first and second order extrapolation the ratio of the errors are exactly equal to the first two dominant Eigen values of the iteration matrix (i.e., 75.7966 and 11.7041). However, the third and higher order extrapolations decompose slowly to successively lower ratios enabling convergence of the Aitken extrapolation.

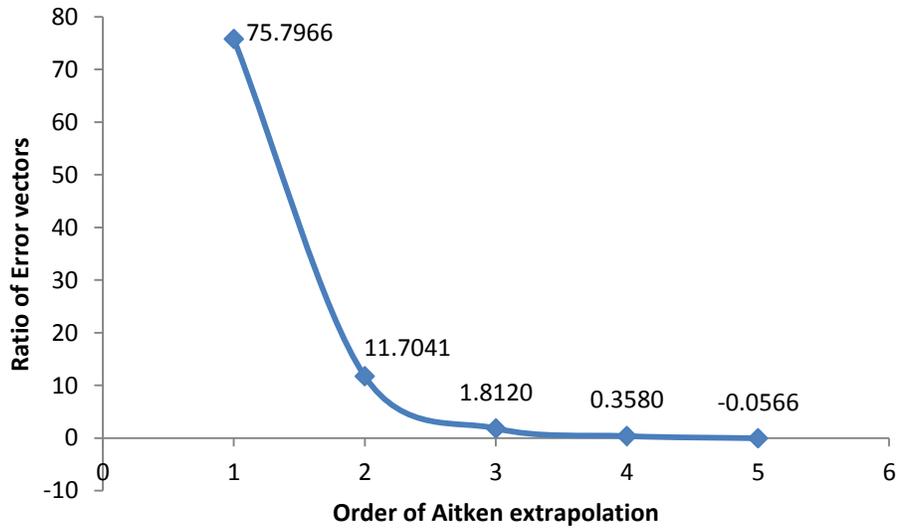

Figure 7. Variation of ratios of the error vectors at 1st to 5th order Aitken extrapolation

The variation of the norm of the error vector with the number of fifth order Aitken iteration is given in Figure 8 below. As can be seen from the figure, the norm of the error vector reduces quickly with the first few iterations.

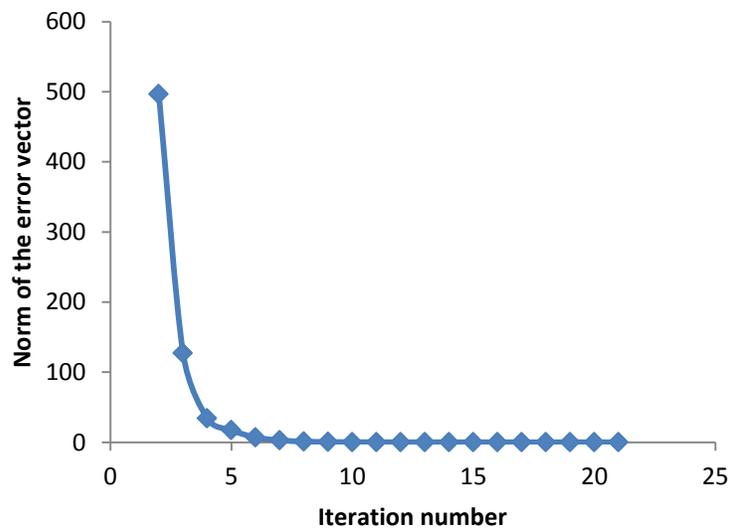

Figure 8. Progress of reduction of the error vector for a 4th order Aitken extrapolation

# 6. Conclusion

The Gauss-Seidel iteration in the case of a convergent iteration is known to follow a diminishing geometric series making it suitable for extrapolation through examination of the ratios of consecutive differences in the solution vector at each step of the iteration. The ratio belongs to the largest Eigen value of the iteration matrix. When sufficient digits of accuracy are obtained for the ratio, the process can be extrapolated towards the solution using a geometric series sum. This procedure is the equivalent of Aitken extrapolation for a convergent iteration. A significant reduction in the number of iteration required is obtained through such extrapolation. Coupling of the successive over relaxation technique with Aitken's extrapolation is possible with further reduction in iteration while employing relaxation factors not necessarily restricted the optimum value which may be difficult to predict in advance for some types of equations.

Coupling of extrapolation with SOR technique is normally not always possible at the optimum acceleration factor w because the largest Eigen value at this optimum value can turn out to be a complex number. Therefore, coupling with SOR technique is done at $\omega$ value typically less than the optimum as the heat flow example presented in this paper showed.

In the case of a divergent Gauss-Seidel iterations the application of Aitken extrapolation formula is made possible and in many cases the extrapolation at each level of the Gauss-Seidel iteration indicates convergence towards the solution. Higher order Aitken extrapolation successively decomposes the dominant Eigen values of the iteration matrix. In doing so, the iteration is successively transformed from an expanding (divergent) form to a stable convergent iteration. At each stage of the application of higher order Aitken extrapolation, the ratio of the error vectors (differences in successive x values) approaches the dominant Eigen value for that order of extrapolation. Higher order Aitken extrapolation provides an interesting possibility of stabilizing, i.e., converting a divergent fixed point iteration to a stable, convergent iteration. In general, the ability to extrapolate the solution from divergent Gauss-Seidel iteration is an interesting possibility that expands further the scope of application of fixed point iteration which in some cases is hampered by problems of divergence.